 \documentclass[11pt]{article}
 \usepackage[top=1in,bottom=1in,left=1.5in,right=1.5in]{geometry}
  \usepackage{amsmath,amssymb}
  \usepackage{latexsym}
 \usepackage[dvips]{pstricks} 
  \usepackage{pst-node}
  \usepackage[dvips]{graphicx}

 \newcommand\R{\mathord{\mathbb R}}
 
 \newcommand\C{\mathord{\mathbb C}}

 \newcommand\F{\mathord{\mathbb F}}

 \newcommand\Z{\mathord{\mathbb Z}}

 \newcommand\N{\mathord{\mathbb N}}

  \renewcommand{\a}{\mathbf{a}}
  \renewcommand{\b}{\mathbf{b}}

  \renewcommand{\c}{\mathbf{c}}
  \renewcommand{\d}{\mathbf{d}}
  
  \newcommand{\e}{\mathbf{e}}
  \newcommand{\f}{\mathbf{f}}
  \newcommand{\g}{\mathbf{g}}
  
  \newcommand{\gl}{\mathbf{GL}\mathnormal}

  \newcommand{\U}{\mathbf{U}}
  \renewcommand{\u}{\mathbf{u}}
  \renewcommand{\v}{\mathbf{v}}
  
  \newcommand{\w}{\mathbf{w}}
  \newcommand{\x}{\mathbf{x}}
  
  \newcommand{\y}{\mathbf{y}}
  
  \newcommand{\z}{\mathbf{z}}
  \newcommand{\0}{\mathbf{0}}
  \newcommand{\1}{\mathbf{1}}

 \renewcommand\f{{\bf f}\mathnormal}

 \newcommand\cF{{\cal F}}
 \newcommand\cG{{\cal G}}

 \newcommand\cP{{\cal P}}
 \newcommand\cQ{{\cal Q}}
 
 \newcommand\cS{{\cal S}}
 \newcommand\cT{{\cal T}}

 \newcommand\rS{{\rm S}}

  \newcommand{\lan}{\langle}
  \newcommand{\ran}{\rangle}
  \newcommand{\an}[1]{\lan#1\ran}
  
  \newcommand{\hs}{\hspace*{\parindent}}
  \newcommand{\proof}{\hs \textbf{Proof.\ }}

  \newcommand{\Gr}{\mathop{\mathrm{Gr}}\nolimits}
  
  \newcommand{\trans}{^\top}
  \newcommand{\qed}{\hspace*{\fill} $\Box$\\}

  \renewcommand{\rS}{\mathrm{S}}

  \newcommand{\bet}{\boldsymbol{\beta}}

  \newcommand{\rank}{\mathrm{rank\;}}
  \newcommand{\brank}{\mathrm{brank}}
  \newcommand{\srank}{\mathrm{srank\;}}
  \newcommand{\sbrank}{\mathrm{sbrank}}
  \newcommand{\Krank}{\mathrm{Krank\;}}

  \newtheorem{theo}{\bfseries \hs Theorem}[section]
  \newtheorem{defn}[theo]{\bfseries \hs Definition}

  \newtheorem{lemma}[theo]{\bfseries \hs Lemma}
  \newtheorem{corol}[theo]{\bfseries \hs Corollary}
  
  \newtheorem{con}[theo]{\bfseries \hs Conjecture}

  \numberwithin{equation}{section} 

 \renewcommand{\span}{\mathrm{span}}

 \renewcommand\dim{{\rm dim\;}}

 \newcommand\Range{{\rm Range\;}}

\begin{document}

 \title{Remarks on the symmetric rank of symmetric tensors}
 \author{
  Shmuel Friedland\footnotemark[1] 
 }
 \renewcommand{\thefootnote}{\fnsymbol{footnote}}

 \footnotetext[1]{
 Department of Mathematics, Statistics and Computer Science,
 University of Illinois at Chicago, Chicago, Illinois 60607-7045,
 USA, \texttt{friedlan@uic.edu}.  This work was supported by NSF grant DMS-1216393.  
 }

 \renewcommand{\thefootnote}{\arabic{footnote}}
 \date{January 21, 2016 }
 \maketitle
 \begin{abstract}
 We give sufficient conditions on a symmetric tensor $\cS\in \rS^d\F^n$ to satisfy the equality: the symmetric rank of $\cS$, denoted as $\srank \cS$, 
is equal to the rank of $\cS$, denoted as $\rank\cS$.  This is done by considering the rank of the unfolded $\cS$ viewed as a matrix $A(\cS)$.  
The condition is: $\rank\cS\in\{\rank A(\cS),\rank A(\cS)+1\}$.  In particular, $\srank\cS=\rank\cS$ for $\cS\in\rS^d\C^n$ for the cases $(d,n)\in\{(3,2),(4,2),(3,3)\}$.
We discuss the analogs of the above results for border rank and best approximations of symmetric tensors.

 \end{abstract}



\section{Introduction}
For a field $\F$ let $\otimes^d\F^n\supset\rS^d\F^n$ denote $d$-mode tensors and the subspace of symmetric tensors on $\F^n$.   Let $\cT\in
\otimes^d\F^n$.   Denote by $\rank \cT$ the rank of the tensor $\cT$.  That is, for $\cT\ne 0$ $\rank \cT$ is the minimal number $k$ such that $\cT$ is a sum of $k$ rank
one tensors. ($\rank 0=0$.)  We say that $\cT$ has a unique decomposition as a sum $\rank\cT$ rank one tensors if this decomposition is unique up to a permutations of the summands.
Assume that $\cS\in\rS^d\F^n\setminus\{0\}$.  Suppose that $|\F|\ge d$, i.e. $\F$ has at least $d$ elements.   Then it is known that $\cS$ is a sum of $k$ symmetric 
rank one tensors \cite[Proposition 7.2]{FS13}.    See \cite{AH95} for the case $|\F|=\infty$, i.e.  $\F$ has an infinite number of elements.  
The minimal $k$ is the symmetric rank of $\cS$, denoted as $\srank \cS$.  
Clearly, $\rank\cS\le \srank\cS$.  In what follows we assume that $d\ge 3$ unless stated otherwise.  
In \cite[P15, page 5]{Oed08} P. Comon asked if $\rank\cS= \srank \cS$ over $\F=\R,\C$.   This problem is also raised in \cite[end \S4.1, p' 1263]{CGLM}. This problem is sometimes referred as Comon's conjecture.  
In \cite{CGLM} it is shown that this conjecture holds in the first nontrivial case: $\rank \cS=2$.

For a finite field the situation is more complicated:  Observe first that for $\F=\Z_2$ and the  symmetric matrix $A=\left[\begin{array}{cc}0&1\\1&0\end{array}\right]$ we have the 
the inequality $\rank A=2 <\srank A=3$.  ($A$ is a sum of all three distinct symmetric rank one matrices in $\rS^2\Z_2^2$.)  Second, it is shown in \cite[Proposition 7.1]{FS13} 
that over a finite field there exist symmetric tensors that are not a sum symmetric rank one tensors.

To state our result we need the following notions:  For $n\in\N$ denote $[n]=\{1,\ldots,n\}$.
Let $\cS=[s_{i_1,\ldots,i_d}]_{i_1,\ldots,i_d\in[n]}\in \rS^d\F^n$.  Denote by $A(\cS)$ an $n\times n^{d-1}$ matrix with entries $b_{\alpha\bet}$ where $\alpha\in[n]$ 
and $\bet=(\beta_1,\ldots,\beta_{d-1})\in [n]^{d-1}$. 
Then $b_{\alpha \bet}:=s_{\alpha,\beta_1,\ldots,\beta_{d-1}}$.  
($A(\cS)\in \F^{n\times n^{d-1}}$ is the unfolding of $\cS$ in the direction $1$.  As $\cS$ is symmetric, the unfolding in every direction $k\in[d]$ gives rise to the same matrix.)  
Hence $\rank A(\cS)\le n$.
If $m:=\rank A(\cS)<n$ it means that we can choose another basis so that $\cS$ is represented as $\cS'\in \rS^d\F^{m}$.
Recall that $\rank \cS\ge \rank A(\cS)$.  (See for example the arguments in \cite{Fri12} for $d=3$.).  
Thus, to study Comon's conjecture we can assume without loss of generality that $\rank A(\cS)=n$.

Denote by $\Sigma(n,d,\F)$ and $\Sigma_s(d,n,\F)$  the Segre variety of rank one tensors plus the zero tensor
and the subvariety of symmetric tensors of at most rank one in $(\F^n)^{\otimes d}$.
 
Let $F_{d,n,k}:\Sigma(n,d,\F)^k\to (\F^n)^{\otimes d}$ be the polynomial map:
\begin{equation}\label{ksegremap}
F_{d,n,k}((\cT_1,\ldots,\cT_k)):=\sum_{j=1}^k \cT_j.
\end{equation}
Let $\cT=F_{d,n,k}((\cT_1,\ldots,\cT_k))$.   
In what follows we say that the decomposition $\cT=\sum_{j=1}^k \cT_j$  is unique if $\rank \cT=k$ and any decomposition of $\cT$ to a sum of $r$ rank one tensors is 
obtained by permuting the order of the summands in $\cT=\sum_{j=1}^k \cT_j$.

Denote by $G_{d,n,k}$ the restriction of the map $F_{d,n,k}$ to $:\Sigma_s(n,d,\F)^k$.
Thus $F_{d,n,k}(\Sigma(n,d,\F)^k)$ and $G_{d,n,k}(\Sigma_s(n,d,\F)^k)$ are the sets of 
of $d$-mode tensors on $\F^n$ tensors of at most rank $k$ and  of symmetric tensors of at most symmetric rank $k$.

Chevalley's theorem yields that $F_{d,n,k}(\Sigma(n,d,\C)^k)$ and $G_{d,n,k}(\Sigma_s(n,d,\C)^k)$ are constructible sets.
Hence the dimension of $G_{d,n,k}(\Sigma_s(n,d,\C)^k)$ is the maximal rank of the Jacobian of the map $G_{d,n,k}$.

$\cS\in\rS^d\C^n$ is said to have a generic symmetric rank $k$ if the following conditions hold:
First, the dimension of the constructible set  $G_{d,n,k}( \Sigma_s(n,d,\C)^k)$ is greater than the dimension of $G_{d,n,k-1}( \Sigma_s(n,d,\C)^{k-1})$.
Second, there exists a strict subvariety $O\subset \Sigma(n,d,\C)^k$,  such that 
$\cS\in G_{d,n,k}( \Sigma_s(n,d,\C)^k\setminus O)$.  Let 
\begin{equation}\label{defknd}
k_{n,d}:=\frac{{n+d-1\choose d}}{n}.
\end{equation}
 
 Chiantini, Ottaviani and  Vannieuwenhoven showed recently \cite{COV15}  that if $\cS\in \rS^d\C^n$ has a generic symmetric rank $k<k_{n,d}$
then $k=\rank S$.  It is much easier to establish this kind of result for smaller values of $k$ using Kruskal's theorem.  See \cite[Theorem 7.6]{FS13}.
 
The aim of this paper is to establish a much weaker result on Comon's conjecture, which does not use the term \emph{generic}.
In particular we show that Comon's conjecture holds for symmetric tensors of at most rank $3$  and for $3$-symmetric tensors of at most rank $5$ over $\C$ .

Our main result is
\begin{theo}\label{maintheo}  Let $d\ge 3$, $|\F|\ge 3$ and $\cS\in\rS^d\F^n$.  Suppose that $\rank \cS\le \rank A(\cS)+1$.
Then $\srank \cS=\rank\cS$.
\end{theo}

We now summarize briefly the content of this paper.
In \S2 we recall Kruskal's theorem on the rank of $3$-tensor.  In \S3 we prove Theorem \ref{maintheo} for the case $\rank \cS=\rank A(\cS)$.  In \S4 we show 
that each $\cS\in\rS^3\F^2$, where $|\F|\ge 3$, satisfies $\srank \cS=\rank\cS$.  In \S5 we prove Theorem \ref{maintheo} in the case $d=3$ and 
$\rank \cS=\rank A(\cS)+1$.   In \S6 we prove Theorem \ref{maintheo} for $d\ge 4$.  In \S7 we summarize our results for $\F=\C$.
In \S8 we discuss two other closely related conjectures:  The first one conjectures that it is possible to replace in Comon's conjecture the ranks with border ranks.
We show that this is true if the border rank of $\cS$ is two. 
The second one conjectures that a best $k$-approximation of symmetric tensor can be chosen symmetric.  For $k=1$, i.e. best rank one approximation, this conjecture holds
and it is a consequence of Banach's theorem \cite{Ban38}.
\section{Kruskal's theorem}
We recall Kruskal's theorem for $3$-tensors and any field $\F$.  For $p$ vectors $\x_1,\ldots,\x_p\in \F^n$ denote by $[\x_1\;\x_2\;\cdots\;\x_p]$ the $n\times p$ matrix whose columns
are $\x_1,\ldots,\x_p$.   Kruskal's rank of $[\x_1\;\x_2\;\cdots\;\x_p]$, denoted as $\Krank(\x_1,\ldots,\x_p)$ is the maximal $k$ such that any $k$ vectors in the set 
$\{\x_1,\ldots,\x_p\}$ are linearly independent.  (If $\x_i=0$ for some $i\in[p]$ then $\Krank(\x_1,\ldots,\x_p)=-\infty$.)
\begin{theo}\label{Kruskalthm} (Kruskal)  Let $\F$ be a field, $r\in\N$ and $\x_i\in \F^m,\y_i\in \F^n,\z_i\in\F^p$ for $i\in[r]$.
Assume that 
\begin{equation}\label{3Tdef}
\cT=\sum_{i=1}^r \x_i\otimes \y_i\otimes \z_i.
\end{equation}
Suppose that
\begin{equation}\label{Kruskalcon}
2r+2\le\Krank(\x_1,\ldots,\x_r)+\Krank(\y_1,\ldots,\y_r)+\Krank(\z_1,\ldots,\z_r).
\end{equation}
Then $\rank \cT=r$.  Furthermore, the decomposition \eqref{3Tdef} is unique.
\end{theo}
Note that $\max(\Krank(\x_1,\ldots,\x_r),\Krank(\y_1,\ldots,\y_r),\Krank(\z_1,\ldots,\z_r))\le r$.
Hence \eqref{Kruskalcon} yields that 
\begin{equation}\label{conkruscon}
\min(\Krank(\x_1,\ldots,\x_r),\Krank(\y_1,\ldots,\y_r),\Krank(\z_1,\ldots,\z_r))\ge 2.
\end{equation}
In particular, $\min(m,n,p)\ge 2$.

In what follows we need a following simple corollary of Kruskal's theorem:
\begin{lemma}\label{Kruslem}  Let $3\le d\in\N$.  Assume that $\x_{j,1},\ldots,\x_{j,r}\in\F^{n_j}$ are linearly independent for each $j\in[d]$.  
Let 
\begin{equation}\label{dTdef}
\cT=\sum_{i=1}^r \otimes_{j=1}^d \x_{j,i}.
\end{equation}
Then $\rank \cT=r$.  Furthermore, the decomposition \eqref{dTdef} is unique.
\end{lemma}
\proof  Observe first that $\otimes_{j=1}^p \x_{j,1},\ldots,\otimes_{j=1}^p \x_{j,r}$ linearly independent for $p=1,\ldots,d$.  Clearly, this is true for $p=1$ and $p=2$.
Use the induction to prove this statement for $p\ge 3$ by observing that  $\otimes_{j=1}^p \x_{j,i}= (\otimes_{j=1}^{p-1} \x_{j,i})\otimes \x_{p,i}$ for $p=3,\ldots,d$.

Consider $\cT$ given by \eqref{dTdef}.   Suppose first that $r=1$.  Then $\cT$ is a rank one tensor and its decomposition is unique.
Assume that $r\ge 2$.
Consider $\cT$ as a $3$-tensor on the $3$-tensor product $\F^{n_1}\otimes \F^{n_2}\otimes(\otimes_{j=3}^d\F^{n_j})$.  
Clearly
\begin{eqnarray}\label{Kequal}
&&\Krank(\x_{1,1},\ldots,\x_{r,1})=\Krank(\x_{1,2},\ldots,\x_{r,2})=\\
&&\Krank(\otimes_{j=3}^d\x_{j,1},\ldots,\otimes_{j=3}^d\x_{j,r})=r.\notag
\end{eqnarray}
As $3r-2\ge 2r$, Kruskal's theorem yields that the rank of $\cT$ as $3$-tensor is $r$.  Hence $\rank \cT$ as $d$ tensor is $r$ too.  
Furthermore the decomposition \eqref{dTdef} of $\cT$ as a $3$-tensor is unique.   Hence the decomposition\eqref{dTdef} is unique.
\qed

In what follows we need the following lemma.
\begin{lemma}\label{rankdecsym}  Let $d\ge 3, n\ge 2$ and $\cS\in\rS^d\F^n$.   Assume that 
\begin{equation}\label{rankdecS}
\cS=\sum_{i=1}^k\otimes_{j=1}^d \x_{j,i}.
\end{equation}  
Then $\cS=\sum_{i=1}^k\otimes_{j=1}^d \x_{\sigma(j),i}$ for any permutation $\sigma$ of $[d]$ .  
Suppose that the following inequality holds:
\begin{eqnarray}\notag
&&2k+2\le K(\x_{1,1},\ldots,\x_{1,k})+ K(\x_{2,1},\ldots,\x_{2,k})+\\ 
&&K(\otimes_{j=3}^d\x_{j,1},\ldots,\otimes_{j=3}^d\x_{j,k}).
\label{kruskuniqsym}
\end{eqnarray}
Then $\rank \cS=\srank\cS=k$, i.e. $\span(\x_{1,i})=\ldots=\span(\x_{d,i})$ for each $i\in[k]$.  Furthermore, the decomposition \eqref{rankdecS}  is unique. 
\end{lemma}
\proof  Assume that \eqref{rankdecS} holds.  Since $\cS$ symmetric we deduce that $\cS=\sum_{i=1}^k\otimes_{j=1}^d \x_{\sigma(j),i}$ for any permutation $\sigma$ of $[d]$ . 
Suppose that \eqref{kruskuniqsym} holds.  Kruskal's theorem yields that the decomposition of $\cS$ as a $3$-tensor on $\F^n\otimes\F^n \otimes (\otimes^{d-2}\F^n)$
is unique.   In particular, the decomposition \eqref{rankdecS} is unique.  Hence $\rank \cS=k$.   Let $\sigma$ be the transposition on $[d]$ satisfying $\sigma(1)=2,\sigma(2)=1$.  
Then $\cS=\sum_{i=1}^k\otimes_{j=1}^d \x_{\sigma(j),i}$.   \eqref{conkruscon} yields that $K(\otimes_{j=3}^d\x_{j,1},\ldots,\otimes_{j=3}^d\x_{j,k})\ge 2$.
That is, the rank one tensors $\otimes_{j=3}^d\x_{j,p}$ and $\otimes_{j=3}^d\x_{j,q}$ are linearly independent for $p <q$.
The uniqueness of the decomposition \eqref{rankdecS},  (up to a permutation of summands), yields that 
$\x_{1,i}\otimes\x_{2,i}\otimes(\otimes_{j=3}^d\x_{j,i}) =\x_{2,i}\otimes\x_{1,i}\otimes(\otimes_{j=3}^d\x_{j,i})$  for each $i\in[n]$.  
Hence $\x_{1,i}\otimes\x_{2,i}=\x_{2,i}\otimes\x_{1,i}$.
Therefore $\x_{1,i}$ and $\x_{2,i}$ are linearly dependent nonzero vectors.  Let $\sigma$ be a transposition on $[d]$ satisfying
$\sigma(2)=p,\sigma(p)=2$ for some $p\ge 3$.   \eqref{conkruscon} yields that $K(\x_{1,1},\ldots,\x_{1,k})\ge 2$.
The uniqueness of the decomposition \eqref{rankdecS}, (up to a permutation of summands), yields that $\x_{1,i}\otimes\x_{2,i}\otimes(\otimes_{j=3}^d\x_{j,i})=\x_{1,i}\otimes\x_{p,i}
\otimes(\otimes_{j=3}^d\x_{\sigma(j),i})$.  Therefore $\x_{2,i}$ and $\x_{p,i}$ are collinear  for each $i\in[n]$. 
Hence  $\span(\x_{1,i})=\ldots=\span(\x_{d,i})$ for each $i\in[d]$.   Thus the decomposition \eqref{rankdecS} is a decomposition to a sum of symmetric rank one tensors.
Hence $\srank\cS=\rank \cS$.\qed

\section{The case  $\rank \cS=\rank A(\cS)$} 
\begin{theo}\label{eqcase}   Let $d\ge 3, n\ge 2$  and $\cS\in\rS(d,\F^n)$.  Suppose that $\rank \cS=\rank A(\cS)$.
Then $\srank \cS=\rank\cS$.  Furthermore, $\cS$ has has a unique rank one decomposition. 
\end{theo} 
\proof  We can assume without loss of generality that $\rank A(\cS)=n$.  So \eqref{rankdecS} holds for $k=n$.   Clearly,  $\x_{j,1},\ldots,\x_{j,n}$ are linearly independent for each $j\in[d]$.  Hence $K(\x_{j,1},\ldots,\x_{j,n})=n$ for $j\in[d]$.
The proof of Lemma  \ref{Kruslem} yields that $K(\otimes_{j=3}^d\x_{j,1},\ldots,\otimes_{j=3}^d\x_{j,n})=n$.  Therefore equality \eqref{Kequal} holds for $r=n$.  As $n\ge 2$ we deduce \eqref{kruskuniqsym} for $k=n$.
Lemma \ref{rankdecsym} yields the theorem.
\qed
 
The following corollary generalizes \cite[Proposition 5.5]{CGLM} to any field $\F$:
\begin{corol}\label{rank2eq}  Let $\F$ be a field, $\cS\in\rS^d\F^n\setminus\{0\}, d\ge 3$.  Assume that $\rank \cS\le 2$.  Then $\srank \cS=\rank\cS$.
\end{corol}
\proof   Clearly, $\rank A(S)\in\{1,2\}$.  If $\rank A(S)=1$ then $\cS=s\otimes^d \u$.  Hence $\rank \cS=\srank \cS=1$. 
 If $\rank A(\cS)=2$ then $\rank \cS=2$ and we conclude the result from Theorem \ref{eqcase}.\qed
\section{The case $\rS^3\F^2$}
\begin{theo}\label{case32}  Let $\cS\in\rS^3\F^2$.  Assume that $|\F|\ge 3$.  Then $\rank \cS=\srank \cS\le 3$.
\end{theo}
For $\F=\C$ this result follows from the classical description of binary forms in two variables due Sylvester \cite{Syl04}.  More generally, consult with \cite{CS11} for results on the rank of tensors in $\rS^d\F^2$ for an algebraic closed field $\F$ of characteristic zero.\\ 

\proof  In view of Corollary \ref{rank2eq} it is enough to consider the case where $\rank \cS\ge 3$.
Let $\cS=[s_{i,j,k}]_{i,j,k\in[2]}$.  
\begin{enumerate}
\item
Assume that $s_{1,1,2} s_{1,2,2}\ne 0$.  Let 
\begin{equation}\label{symS32dec}
\cS=\sum_{i=1}^3 t_i\otimes^3\u_i, \;\u_1=(1,b)\trans,\quad \u_2=(1,0)\trans,\quad \u_3=(0,1)\trans.
\end{equation}
Then 
\begin{eqnarray*}&s_{1,1,2}=t_1 b, \quad s_{1,2,2}=t_1 b^2 \Rightarrow  t_1=\frac{s_{1,1,2}^2}{s_{1,2,2}},\quad  b=\frac{s_{1,2,2}}{s_{1,1,2}},\\
&t_2=s_{1,1,1}-t_1, \quad t_3=s_{2,2,2}-t_1b^3.
\end{eqnarray*}
Hence $\rank\cS\le 3$.  Our assumption yields that $\rank \cS=3$ and \eqref{symS32dec} is a minimal decomposition of $\cS$ to rank one tensors.
This decomposition shows that $\rank \cS=\srank \cS$.
\item  Assume that  $s_{1,1,2}=s_{1,2,2}=0$ then $\cS=s_{1,1,1}\otimes^3 (1,0)\trans
+s_{2,2,2}\otimes^3(0,1)\trans$.  This contradicts our assumption that $\rank \cS\ge 3$.
\item
It is left to discuss the case where $\rank \cS\ge 3$ and  $s_{1,1,2}=0$ and $s_{1,2,2}\ne 0$.  The homogeneous polynomial of degree $3$ corresponding to $\cS$
is 
\[f(x_1,x_2)=s_{1,1,1}x_1^3 +s_{1,2,2}x_1x_2^2 +s_{2,2,2}x_2^3.\]
\begin{enumerate}
\item
Assume that  the characteristic of $\F$ is $3$.  Make the following change of variables: $x_1=y_1, x_2=y_1+y_2$.  The new tensor $\cS'$ satisfies 
\emph{1}.
\item Assume that  the characteristic of $\F$ is not $3$.  
\begin{enumerate}
\item Assume that $s_{1,1,1}\ne 0$.  
Make the following change of variables:  $x_1=y_1+ay_2, x_2=y_2$.  Then
\[f(y_1,y_2)=\alpha y_1^3 +\beta y_1^2 y_2 +\gamma y_1y_2^2 +\delta, \; \beta=3a s_{1,1,1},\gamma=s_{1,2,2}+3a^2s_{1,1,1}.\]
Then choose a nonzero $a$ such that $s_{1,2,2}+3a^2s_{1,1,1}\ne 0$.   (This is always possible if $|\F|\ge 4$ as we assumed that $\F\ne \Z_3$ and $|\F|\ge 3$.) 
The new tensor $\cS'$ satisfies 
\emph{1}.
\item Assume that $s_{1,1,1}=s_{2,2,2}=0$. Make the following change of variables: $x_1=y_1, x_2=(y_1+y_2)$. Then we are
either in the case \emph{1} if the characteristic of $\F$ is not $2$ or in the case \emph{3}(b)i if the characteristic of $\F$ is $2$.  
\item   Assume that $s_{1,1,1}=0$ and $s_{2,2,2}\ne 0$.  
Make the following change of variables: $y_2=s_{1,2,2}x_1+s_{2,2,2}x_2,y_2=x_2$.   Then we are in the case \emph{3}(b)ii.\qed
\end{enumerate}
\end{enumerate}
\end{enumerate}

Note that if $|\F|\gg 1$ then using the change of coordinates and then the above procedure we obtain that if $\rank \cS=3, \rank A(\cS)=2$ we have many presentation of
$\cS$ as sum of three rank one symmetric tensors.  

Observe next that for $\F=\Z_2$ not every symmetric tensor $\cS\in\rS^3\Z_2^2$ is a sum of rank one symmetric tensors.
The number of all symmetric tensors in $\rS^3\Z_2^2$ is $2^4$.  The number of all nonzero symmetric tensors which are sum of rank one symmetric tensors is $2^3-1$.
Hence Theorem \ref{case32} does not hold for $\F=\Z_2$.
\begin{corol}\label{rank3symten}  Let $\cS\in\rS^3\F^n$.  Assume that $|\F|\ge 3$ and $\rank \cS=3$.  Then $\srank\cS=\rank \cS$.
\end{corol}
\proof  Clearly, $\rank A(\cS)\in\{2,3\}$.  If $\rank A(\cS)=3$ we deduce the corollary from Theorem \ref{eqcase}.  If $\rank A(\cS)=2$ we deduce the corollary from Theorem
\ref{case32}.\qed
\section{The case $d=3$ and $\rank \cS=\rank A(\cS)+1$}
In this section we prove Theorem \ref{maintheo} for $d=3$.  In view of Theorem \ref{eqcase} it is enough to consider the case $\rank \cS=\rank A(\cS)+1$.
Furthermore, in view of Theorem \ref{case32} it is enough to consider the case $\rank A(\cS)\ge 3$.
 We first give the following obvious lemma:
 \begin{lemma}\label{3tenthlemplus1}   Let $|\F|\ge 3$ and $\cS\in \rS(3,\F^n)$.  Suppose that $\rank \cS= \rank A(\cS)+1$.
 Assume furthermore that there exists a decomposition of $\cS$ to $\rank A(\cS)+1$ rank one tensors such that at least one of them is symmetric, i.e.
 $s\otimes^3 \u$.  Let $\cS'=\cS-s\otimes^3\u$.  Then $\rank \cS'=\rank \cS-1$ and $\rank A(\cS')\in\{ \rank A(\cS)-1, \rank A(\cS)\}$.   Furthermore:
 \begin{enumerate}
\item  If $\rank A(\cS')=\rank A(\cS)$ then $\rank \cS'=\rank A(\cS')$
and $\rank \cS'=\srank\cS'$.  Hence $\rank \cS=\srank\cS$.
\item If $\rank A(\cS')= \rank A(\cS)-1$ then $\rank \cS'=\rank A(\cS')+1$.
\end{enumerate}
\end{lemma}

\subsection{The case $\rank A(\cS)=3$}
We now discuss Theorem \ref{maintheo}  where $\cS\in\rS(3,\F^n)$, where $\rank\cS=4, \rank A(\cS)=3$.  Without loss of generality we can assume that $n=3$.  Then
\begin{equation}\label{rankdecS34}
\cS=\sum_{i=1}^4 \x_i\otimes\y_i\otimes \z_i.
\end{equation}

Suppose first that there is a decomposition \eqref{rankdecS34} such that  $\x_i\otimes\y_i\otimes \z_i$ is symmetric for some $i\in[4]$.  Then we can use Lemma \ref{3tenthlemplus1}.  Apply Theorems \ref{eqcase} and \ref{case32} to deduce that $\srank\cS=\rank\cS=4$.

Assume the \emph{Assumption}: there no is a decomposition \eqref{rankdecS34} such that  $\x_i\otimes\y_i\otimes \z_i$ is symmetric for some $i\in[4]$.   
The first part of Lemma \ref{rankdecsym} yields:
\[0=\sum_{i=1}^4 \x_i\otimes (\y_i\otimes \z_i-\z_i\otimes \y_i).\]
As $\rank A(\cS)=3$ we can assume that $\x_1,\x_2,\x_3$ are linearly independent. 
 Then $\x_4=\sum_{j=1}^3 a_j\x_j$.  
Then the above equality yields:
\begin{equation}\label{4tenseq}
0=\sum_{j=1}^3\x_j\otimes(\y_j\otimes \z_j-\z_j\otimes \y_j +a_j(\y_4\otimes \z_4-\z_4\otimes \y_4)).
\end{equation}
As $\x_1,\x_2,\x_3$ are linearly independent it follows that 
\begin{equation}\label{skewsymeq}
\y_j\otimes \z_j-\z_j\otimes\y_j=-a_j(\y_4\otimes \z_4-\z_4\otimes \y_4) \textrm{ for } j\in[3].
\end{equation}

Assume first that $\y_4$ and $ \z_4$ are collinear.    
Then $\y_j$ and $\z_j$ are collinear for $j\in[3]$.  
Hence we w.l.o.g we can assume that $\y_i=\z_i=\u_i$ for $i\in[4]$.  So $\cS=\sum_{i=1}^4 \x_i\otimes \u_i\otimes \u_i$.  Since $\cS$ is symmetric we can
we obtain that $\cS=\sum_{i=1}^4 \u_i\otimes \u_i\otimes \x_i$.  Renaming the vectors we can assume that in the original decomposition \eqref{rankdecS34} we 
have that $\x_i$ and $\y_i$ are collinear for $i\in [4]$.  Since we assumed that no rank one tensor $\x_i\otimes\y_i\otimes\z_i$ is not symmetric, we deduce that 
each pair $\y_i,\z_i$ in the original decomposition \eqref{rankdecS34} is not collinear. In particular, it is enough to study the case where  $\y_4$ and $\z_4$ are not collinear, i.e. $\y_4\otimes \z_4-\z_4\otimes \y_4\ne 0$.  Let $\U:=\span(\y_4,\z_4)$.  Then $\dim\U=2$.

Suppose that $a_j\ne 0$ for some $j\in[3]$.
Then \eqref{skewsymeq} yields that 
$\y_j$ and $\z_j$ are not collinear and $\span(\y_j,\z_j)= \U$.  Assume that $a_j=0$.  Then $\y_j$ and $\z_j$ are collinear. 

Assume first that $a_1a_2a_3\ne 0$.  Then the above arguments yields that
$\span(\y_1,\ldots,\y_4)$ $\subseteq \U$, which contradicts the assumption that $\span(\y_1,\ldots,\y_4)=\F^3$.

So we need to assume that at least one of $a_i=0$.
Assume first that exactly one $a_i=0$.  Without loss of generality we can assume that    
in \eqref{4tenseq} $a_1=0$ and $a_2 a_3\ne 0$.
This yields that $\y_1$ and $\z_1$ are collinear.  
Our \emph{Assumption} yields that  $\x_1$ and $\y_1$ are not collinear.
Furthermore $\span(\y_2,\z_2)=\span(\y_3,\z_3)=\U$.   Hence $\span(\y_2,\y_3,\y_4)\subseteq \U$.
As span $\y_1,\y_2,\y_3,\y_4$ is the whole space, 
we deduce that $\y_1\not\in\span(\y_2,\y_3,\y_4)$.  Similarly $\y_1\not\in\span(\z_2,\z_3,\z_4)$.
Furthermore, $\dim\span(\y_2,\y_3,\y_4)=2$.  Hence $\span(\y_2,\y_3,\y_4)=\U$.
We now recall that $\cS=\sum_{i=1}^4 \y_i\otimes\x_i\otimes \z_i$.  
Again, by renaming the indices $2,3,4$ we can assume that $\y_1,\y_2,\y_3,$ are linearly independent.  Since $\span(\y_2,\y_3,\y_4)=\U$ it follows that $\y_4=b_2\y_2+b_3\y_3$.  Since $\cS$ is symmetric we have the equality $\cS=\sum_{i=1}^4 \y_i\otimes \x_i \otimes \z_i$.  Permuting the last two factors
we obtain the equality $0=\sum_{i=1}^4 \y_i\otimes (\x_i\otimes \z_i-\z_i\otimes \y_i)$.  Hence we have an analogous equality to \eqref{4tenseq}: 
\[0=\y_1\otimes (\x_1\otimes \z_1-\z_1\otimes \x_1)+\sum_{j=2}^3 \y_j\otimes(\x_j\otimes \z_j-\z_j\otimes \x_j +b_j(\x_4\otimes \z_4-\z_4\otimes \x_4)).\]
Therefore $\x_1\otimes \z_1-\z_1\otimes \x_1=0$.   Thus $\x_1$ and $\z_1$ are collinear.  Recall that we already showed that $\y_1$ and $\z_1$ are collinear.  Hence $\x_1\otimes\y_1\otimes \z_1$ is a symmetric rank one tensor.  
 So we have a contradiction to our \emph{Assumption}.
 
 Finally let us assume that $a_i=a_j=0$ for some two distinct indices $i,j\in[3]$.  W.l.o.g. we can assume that $\x_4=\x_3$, i.e. $a_1=a_2=0, a_3=1$.  
This implies that $\y_i$ and $\z_i$ are collinear for $i=1,2$.  Furthermore 
\[C:=\y_3\otimes \z_3+\y_4\otimes \z_4=\z_3\otimes \y_3+ \z_4\otimes \y_4.\]
So $C$ is  a symmetric matrix.  Note that $C$ is a rank two matrix. Otherwise $\y_3\otimes \z_3$ and $\y_4\otimes \z_4$ are collinear.
Then $\x_3\otimes \y_3\otimes \z_3+\x_3\otimes \y_4\otimes \z_4$ is a rank one tensor.  So $\rank \cS\le 3$, contrary to our assumptions.
Thus we can assume that $C=\y_3\otimes \y_4+\y_4\otimes \y_3$ and $\y_3,\y_4$ are linearly independent.  Hence we can assume that
\begin{eqnarray*}
\cS=&&\x_1\otimes\y_1\otimes\y_1+\x_2\otimes\y_2\otimes\y_2+\x_3\otimes(\y_3\otimes\y_4+\y_4\otimes\y_3)=\\
&&\y_1\otimes\x_1\otimes\y_1+\y_2\otimes\x_2\otimes\y_2+\y_3\otimes\x_3\otimes\y_4+\y_4\otimes \x_3\otimes\y_3=\\
&&\y_1\otimes\y_1\otimes\x_1+\y_2\otimes\y_2\otimes\x_2+\y_3\otimes\y_4\otimes\x_3+\y_4\otimes \y_3\otimes\x_3.
\end{eqnarray*}
Our \emph{Assumption} yields that the pairs $\x_1,\y_1$ and $\x_2,\y_2$ are linearly independent. Hence $Q:=\x_1\otimes\y_1-\y_1\otimes \x_1\ne 0$.
Subtracting the third expression for $\cS$ from the second one we deduce 
\begin{eqnarray*}
&&\y_1\otimes (\x_1\otimes\y_1-\y_1\otimes\x_1)+\y_2\otimes (\x_2\otimes\y_2-\y_2\otimes\x_2)+\\
&&\y_3\otimes (\x_3\otimes\y_4-\y_4\otimes\x_3)+
\y_4\otimes (\x_3\otimes\y_3-\y_3\otimes\x_3)=0
\end{eqnarray*}
As  $\y_3,\y_4$ are linearly independent, without loss in generality we may assume that $\y_2,\y_3,\y_4$ are linearly independent.  So $\y_1=b_2\y_2+b_3\y_3+b_4\y_4$.  Substitute in the above equality this expression for $\y_1$ only for the $\y_1$ appearing in the left-hand side to obtain
\[\y_2\otimes(\x_2\otimes \y_2-\y_2\otimes\x_2+b_2 Q_2)+\y_3(\x_3\otimes \y_4-\y_4\otimes \x_3+b_3 Q)+\y_4(\x_3\otimes\y_3-\y_3\otimes \x_3+b_4Q)=0.
\]
Hence
\[\x_2\otimes \y_2-\y_2\otimes\x_2+b_2 Q_2=\x_3\otimes \y_4-\y_4\otimes \x_3+b_3 Q=\x_3\otimes\y_3-\y_3\otimes \x_3+b_4Q=0.\]
Note that our \emph{Assumption} yields that $b_2\ne 0$.   Hence $\span(\x_2,\y_2)=\span(\x_1,\y_1)$.
Suppose first that $b_3\ne 0$.  Then 
$\span(\x_3,\y_4)=\span(\x_1,\y_1)$.  This contradicts the assumption that $\x_1,\x_2,\x_3$ are linearly independent.  As $\x_3=\x_4$, we get also a contradiction if $b_4\ne 0$.
Hence $b_3=b_4=0$.  So $\y_3,\y_4\in\span(\x_3)$.  This contradicts the assumption that $\y_3$ and $\y_4$ are linearly independent.

In conclusion we showed that our \emph{Assumption} never holds.
The proof of this case of Theorem \ref{maintheo} is concluded.\qed

\subsection{Case $\rank A(\cS)\ge 4$}
\proof  By induction on $r=\rank A(\cS)\ge 3$.    For $r=3$ the proof follows from the results above.  Assume that Theorem holds for $\rank \cS=r+1$.  Assume now that $\rank A(\cS)=r+1$ and $\rank \cS=r+2$.  Without loss of generality we can assume that $n=r+1$.
Suppose first that  the assumptions of Lemma \ref{3tenthlemplus1} hold.
If we are in the case \emph{1} then $\srank \cS=\rank\cS$.  If we are in the case \emph{2}. then we deduce from the induction hypothesis that 
$\srank\cS=\rank\cS$.

As in the proof of the case $\rank\cS=3$ we assume the \emph{Assumption}: 
There does not exist a decomposition of $\cS$ to $\rank A(\cS)+1$ rank one tensors such that at least one of them is symmetric.
We will show that we will obtain a contradiction.

Suppose $\cS=\sum_{i=1}^{n+1} \x_i\otimes\y_i\otimes\z_i$.  
The we have the equality
\[0=\sum_{i=1}^{n+1} \x_i\otimes( \y_i\otimes \z_i-\z_i\otimes \y_i)\]
and the fact that span of all x's, y's and z's is $\F^n$. 

Without loss of generality we may assume that $\x_1,\ldots,\x_{n}$ are linearly independent.
 So $\x_{n+1}=\sum_{i=1}^{n} a_i\x_i$.  Hence
\begin{equation}\label{bas3drel}
\sum_{i+1}^{n}\x_i\otimes(\y_i\otimes \z_i-\z_i\otimes \y_i + a_i(\y_{n+1}\otimes \z_{n+1}-\z_{n+1}\otimes \y_{n+1}))=0.
\end{equation}
As in the case $n=3$ we can assume that $\y_{n+1}$ and $\z_{n+1}$ are not collinear.
Thus if $a_i=0$ we deduce that $\y_i$ and $\z_i$ are collinear.
If $a_i\ne 0$ we deduce that $\span(\y_i,\z_i)=\span(\y_{n+1},\z_{n+1})$.  Since $\y_1,\ldots,\y_{n+1}$ span $\F^n$ we can have at most two nonzero $a_i$.
Since $\x_{n+1}\ne 0$ we must have at least one nonzero $a_i$.  Assume first that $n-1$ out of $\{a_1,\ldots,a_n\}$ are zero.  We may assume without loss of generality
that $a_1=\ldots=a_{n-1}=0$ and $a_n=1$.  So $\x_{n+1}=\x_n$.  Without loss of generality we may assume that
\[\cS=\x_n\otimes (\y_n\otimes\z_n+\y_{n+1}\otimes\z_{n+1})+ \sum_{i=1}^{n-1} \x_i\otimes \y_i\otimes\y_i.\]
Since $\cS$ is symmetric as in case $n=3$ we deduce that $\y_n\otimes\z_n+\y_{n+1}\otimes\z_{n+1}$ is symmetric and has rank two.
So we can assume that $\z_{n}=\y_{n+1}, \z_{n+1}=\y_n$ and $\dim\span(\y_n,\y_{n+1})=2$.  We now repeat the arguments in the proof of this case for $n=3$ to deduce
the contradiction.

Suppose finally that exactly $n-2$ out of $\{a_1,\ldots,a_n\}$ are zero.  We may assume without loss of generality 
that $a_1=\ldots=a_{n-2}=0$ and $a_{n-1},a_n\ne 0$.  So $\span(\y_{n-1},\z_{n-1})=\span(\y_{n},\z_{n})=\span(\y_{n+1},\z_{n+1})$.
Hence $\y_{n-1},\y_n,\y_{n+1}$ are linearly dependent.  Since $\y_1,\ldots,\y_{n+1}$ span the whole space we must have that $\dim\span(\y_{n-1},\y_n,\y_{n+1})=2$.
Without loss of generality we may assume the following: First, $\y_{n},\y_{n+1}$ are linearly independent and $\y_{n-1}=a\y_n+b\y_{n+1}$.  Second $\z_k=\y_k$ for
$k=1,\ldots,n-2$.  So we can assume that
\[\cS=\sum_{j=1}^{n-2} \x_j\otimes \y_j\otimes \y_j+\sum_{j=n-1}^{n+1} \x_j\otimes \y_j\otimes \z_j=\sum_{j=1}^{n-2} \y_j\otimes \x_j\otimes \y_j+\sum_{j=n-1}^{n+1} \y_j\otimes \x_j\otimes \z_j.\] 
Permuting the las two factors in the last part of the above identity we obtain:
\[0=\sum_{j=1}^{n-2}\y_j\otimes(\x_j\otimes\y_j-\y_j\otimes\x_j)+\sum_{j=n-1}^{n+1}
\y_j(\x_j\otimes\z_j-\z_j\otimes\x_j).\]
Substitute $\y_{n-1}=a\y_n+b\y_{n+1}$ and recall that $\y_1,\ldots,\y_{n-2},\y_{n},\y_{n+1}$ are linearly independent.  Hence $\x_i$ and $\y_i$ are collinear for $i=1,\ldots,n-2\ge 2$.  This contradicts our \emph{Assumption}.\qed
\section{Theorem \ref{maintheo} for $d\ge 4$}
In this section we show Theorem for \ref{maintheo} for $d\ge 4$.  Theorem \ref{eqcase} yields that it is enough to consider the case where $\rank\cS=\rank A(\cS)+1$.
We need the following lemma:
\begin{lemma}\label{linindlemm}  Let $2\le d\in\N$.  Assume that $\x_{j,1},\ldots\x_{j,n+1}\in\F^{n}\setminus\{\0\}$ and 

\noindent
$\span(\x_{j,1},\ldots,\x_{j,n+1})
=\F^n$ for $j\in[d]$.
Consider the $n+1$ rank one $d$-tensors $\otimes_{j=1}^d \x_{j,i}, i\in[n+1]$. 
Then either all of them are linearly independent or $n$ of these tensors are linearly independent and the other one is a multiple of one of the $n$ linearly independent tensors.
\end{lemma}
\proof  It is enough to consider the case where  the $n+1$ rank one $d$-tensors $\otimes_{j=1}^d \x_{j,i}, i\in[n+1]$ are linearly dependent.
Without loss of generality we may assume that $\x_{1,1},\ldots,\x_{1,n}$ are linearly independent.  Hence the $n$ tensors $\otimes_{j=1}^d \x_{j,i}, i\in[n]$ are linearly independent as rank one matrices $\x_{1,i}\otimes (\otimes_{j=2}^d\x_{j,i})$ for $i\in[n]$.   (I.e., the corresponding unfolding of $n$ tensors in mode $1$ are linearly independent.)  Assume that $\x_{1,n+1}=\sum_{j=1}^n a_j\x_{1,j}$ where not $a_j$ are zero.  Since we assumed that $\otimes_{j=1}^d \x_{j,i}, i\in[n+1]$ are linearly dependent it follows that $\otimes_{j=1}^d \x_{j,n+1}=\sum_{i=1}^n b_i\otimes_{j=1}^d \x_{j,i}$.  So we obtain the identity $\sum_{i=1}^n \x_{1,i}\otimes 
\cT_i=0$.  Here $\cT_i\in \otimes^{d-1}\F^n$ is a tensor of at most rank 2.
Since $\x_{1,1},\ldots,\x_{1,n}$ are linearly independent if follows that each $\cT_i$ is zero.  Hence if $a_i\ne 0$ it follows that $b_i$ is not zero and $\otimes_{j=2}\x_{j,n+1}$
and $\otimes_{j=2}\x_{j,i}$ are collinear.  Therefore  $\x_{j,i}$ and $\x_{j,n+1}$ are collinear for $j=2,\ldots,d$.  Since $\dim\span(\x_{j,1},\ldots,\x_{j,n+1})=n$, we can't have another $a_k\ne 0$.  So $\otimes_{j=1}^d \x_{j,n+1}$ is collinear with $\otimes_{j=1}\x_{j,i}$ as we claimed.\qed

\textbf{Proof of Theorem \ref{maintheo} for $d\ge 4$ and $\rank \cS=\rank A(\cS)+1$}. 
Without loss of generality we may assume that $n=\rank A(\cS)\ge 2$. 
Assume that $\cS=\sum_{i=1}^{n+1} \otimes_{j=1}^d \x_{j,i}$.  Clearly, the assumptions of Lemma \ref{linindlemm} holds.  Consider the $d-2$ rank one tensors
$\otimes_{j\in[d]\setminus\{p,q\}}\x_{j,i}$ for fixed $p\ne q\in [d]$ and  $i\in[n+1]$.
Suppose that these $n+1$ rank one tensors are linearly independent.  We claim that $\x_{p,i}$ and $\x_{q,i}$ are collinear for each $i\in[n+1]$.  Without loss of generality we 
may assume that $p=1, q=2$.  By interchanging the first two factors in the representation of $\cS$ as a rank $n+1$ tensor we deduce:
\[\sum_{i=1}^{n+1} (\x_{1,i}\otimes \x_{2,i} - \x_{2,i}\otimes\x_{1,i})\otimes(\otimes_{j=3}^d \x_{j,i})=0.\]
As $\otimes_{j=3}^d \x_{j,i}, i\in[n+1]$ are linearly independent we deduce that $\x_{1,i}\otimes \x_{2,i} - \x_{2,i}\otimes\x_{1,i}=0$ for each $i\in [n+1]$.  I.e., 
$\x_{1,i}$ and $\x_{2,i}$ are collinear for each $i\in[n+1]$.  

Suppose first that for each pair of integers $1\le p<q\le d$  $\otimes_{j\in[d]\setminus\{p,q\}}\x_{j,i}, i\in[n+1]$  are linearly independent.
Hence $\x_{j,i}\in\span(\x_{1,i})$ for $j\in[d]$ and $i\in[n+1]$.  Therefore $\otimes_{j=1}^d \x_{j,i}$ is a rank one symmetric tensor for each $i\in[n+1]$.
Thus $\srank S=\rank S$.

Assume now, without loss of generality, that $\otimes_{j=3}^d\x_{j,i}, i\in[n+1]$  are linearly dependent.  By applying Lemma \ref{linindlemm} we can assume without loss of generality that $\otimes_{j=3}^d\x_{j,i}, i\in[n]$ are linearly independent and
$\otimes_{j=3}^d \x_{j,n+1}=\otimes_{j=3}^d \x_{j,n}$.   Without loss of generality we may assume that $\x_{j,n+1}=\x_{j,n}$ for $j\ge 3$.  (We may need to rescale the vectors $\x_{2,n+1},\ldots,\x_{d,n+1}$.)
Hence $\x_{j,1},\ldots,\x_{j,n}$ are linearly independent for each $j\ge 3$.  Therefore we have the following decomposition of $\cS$ as a $3$-tensor in $(\otimes^2\F^{n})\otimes \F^n \otimes (\otimes^{d-3}\F^n)$:
 \begin{eqnarray}\label{Sdec3ten}
\cS=&&
(\x_{1,n}\otimes \x_{2,n}+\x_{1,n+1}\otimes\x_{2,n+1})\otimes \x_{3,n}\otimes(\otimes_{j=4}^d \x_{j,n})+\\
&&\sum_{i=1}^{n-1} (\x_{1,i}\otimes\x_{2,i})\otimes \x_{3,i}\otimes (\otimes_{j=4}^d \x_{j,i})=\sum_{i=1}^n \cT_i\otimes (\otimes_{j=4}^d \x_{j,i}).\notag
\end{eqnarray}
Clearly, $\otimes_{j=4}^d\x_{j,1},\ldots,\otimes_{j=4}^d\x_{j,n}$ are linearly independent.  Since $\cS$ is symmetric by interchanging every two distinct factors $p,q\in[3]$  in $\otimes^d\F^n$ we deduce that $\cT_1,\ldots,\cT_n$ are symmetric $3$-tensors.  Consider the symmetric tensor $\cT_n=(\x_{1,n}\otimes \x_{2,n}+\x_{1,n+1}\otimes\x_{2,n+1})\otimes \x_{3,n}$.  As  the rank of $A(\cT_n)$ in the the third coordinate is $1$ it follows that $\rank A(S)=1$.  Hence $\rank \cT_n=1$.
Therefore $\rank \cS\le n$ contrary to our assumptions.\qed

\begin{corol}\label{rank3case}  Let $|\F|\ge 3, d\ge 3, n\ge 2$.  Assume that $\cS\in\rS^d\F^n$. Then $\srank\cS=\rank \cS$ under the following assumptions:
\begin{enumerate}
\item $\rank \cS \le 3$.
\item $\srank \cS\le 4$.
\end{enumerate}
\end{corol}
\proof  It is enough to consider the case where $n=\rank A(\cS)\ge 2$.  

\noindent
\emph{1.} Clearly, $\rank \cS\in\{2,3\}$.  Theorem \ref{maintheo} yields that $\srank \cS=\rank \cS$.

\noindent
\emph{2.}  Assume to the contrary that $\rank\cS<\srank \cS\le 4$.  Then $\rank \cS\le 3$.  Part \emph{1.} implies the contradiction $\rank\cS=\srank \cS$.
\qed

\section{Symmetric tensors over $\C$}
Recall the known maximal value of the symmetric rank in $\rS^d\C^n$, denoted as $\mu(d,n)$: 
\begin{enumerate}
\item $\mu(d,2)=d$ \cite{CS11}, \cite[\S3.1]{BT14};
\item $\mu(3,3)=5$ \cite [\S96]{Seg42}, \cite{CM96} and \cite{LT10};
\item $\mu(3,4)=7$ \cite[\S97]{Seg42};
\item $\mu(3,5)\le 10$ \cite{DeP15};
\item $\mu(4,3)=7$ \cite{Kle99, DeP13}.
\end{enumerate}
\begin{theo}\label{valcomocn}  Let $\F=\C$ and $\cS$ be a symmetric tensor in $\rS^d\C^n$.  Then $\srank \cS=\rank \cS$ in the following cases: 
\begin{enumerate}
\item $d\ge 3, n\ge 2$ and $\rank \cS\in\{\rank A(\cS),\rank A(\cS)+1\}$.
\item For $n=2$ and $d=3$.
\item For $n=2$ and $d=4$
\item $n=d=3$.
\item $\cS\in\rS^3\C^n$ and $\rank \cS\le 5$.
\item  $\cS\in\rS^3\C^n$ and $\srank \cS\le 6$.

\end{enumerate}
\end{theo}
\proof  Assume that $\cS\in\rS^d\C^n$.  Clearly, it is enough to prove the theorem for the case  $\rank A(\cS)\ge 2$.
Furthermore, we can assume that $n=A(\cS)$.
Thus it is enough to assume the following conditions:
\begin{equation}\label{obvrankin}
2\le n=\rank A(\cS)\le \rank \cS\le \srank \cS\le \mu(d,n).
\end{equation}

\noindent
\emph{1.}  follows from Theorem \ref{maintheo}.

\noindent
\emph{2.}  Assume $\cS\in\rS^3\C^2$.  As $\mu(3,2)=3$ we deduce the theorem from \emph{1.}

\noindent
\emph{3.}  Assume that $\cS\in\rS^4\C^2$.  Suppose that $\rank \cS\in\{2,3\}$.
Then \emph{1.} yields that $\srank \cS=\rank\cS$.  

Suppose that $\rank\cS\ge 4$.  As $\mu(4,2)=4$ in view of \eqref{obvrankin} it follows that $\srank \cS=\rank\cS=4$.

\noindent
\emph{4.}
Assume now that $\cS\in\rS^3\C^3$.  Suppose first that $\rank A(\cS)=2$.  Then by changing a basis in $\C^3$ we can assume that $\cS\in\rS^3\C^2$.  Part \emph{2.} yields that
$\srank\cS=\rank\cS$.

Suppose that $\rank A(\cS)=3$.  If $\rank\cS\in\{3,4\}$ then \emph{1.} yields that $\srank\cS=\rank\cS$.   Suppose now that $\rank\cS\ge 5$.   \eqref{obvrankin} yields that 
$\srank\cS\ge 5$.  The equality $\mu(3,3)=5$ yields that $\rank\cS=5$.   Hence $\srank \cS=\rank\cS=5$.

\noindent  
\emph{5.} \eqref{obvrankin} yields $\rank A(\cS)\le 5$.   If $\rank A(\cS)=2$ then \emph{2.} yields that $\srank\cS=\rank\cS$.
If $\rank A(\cS)=3$ then \emph{4.} yields that $\srank\cS=\rank\cS$. If $\rank A(\cS)\ge 4$ then \emph{1.} yields that $\srank\cS=\rank\cS$. 

\noindent  
\emph{6.}  Assume to the contrary that $\rank\cS<\srank\cS$.  So $\rank\cS\le 5$.  \emph{5.} implies the contradiction $\rank\cS=\srank\cS$.
\qed
 
\section{Two version of Comon's conjecture}\label{S:bc}
In this section we assume that $\F=\R,\C$.
\subsection{Border rank}
\begin{defn}\label{defbordrank}
Let $\cT\in\otimes^d \F^n\setminus\{0\}$.  Then the border of $\cT$, denoted as $\brank_{\F}\;\cT$, is $r \in \N$ if the following conditions hold
\begin{enumerate}
\item There exists a sequence $\cT_k\in\otimes^d\F^n,k\in\N$ such that $\rank \cT_k=r$ for $k\in\N$ and $\lim_{k\to\infty}\cT_k=\cT$.
\item Assume that a sequence $\cT_k\in\otimes^d\F^n,k\in\N$ converges to $\cT$.  Then 

\noindent
$\liminf_{k\to\infty}\rank\cT_k\ge r$.
\end{enumerate}
\end{defn}
Clearly, $\brank_{\F}\;\cT\le \rank \cT$.  For $d=2$ it is well known that $\brank_{\F}\;\cT=\rank\cT$.  
Hence 
\begin{equation}\label{brankATin}
\rank A(\cT)\le \brank_{\F}\;\cT.
\end{equation}
For $d>2$ one has examples where  $\brank_{\F}\;\cT<\rank\cT$
\cite{CGLM}:  Assume that $\x,\y\in\F^n$ are linearly independent.  Let 
\begin{equation}\label{ex3symtesrnk3}
\cS=\x\otimes\x\otimes\y+\x\otimes \y\otimes\x+\y\otimes\x\otimes\x, \quad \cS=\lim_{\epsilon\to 0} \frac{1}{\epsilon} (\otimes^3(\x+\epsilon\y) -\otimes^3 \x).
\end{equation}
It is straightforward to show that $\rank \cS=3,\brank_{\F}\;\cS=2$.   (See the proof of Theorem \ref{brank2}.)

Assume that $\cS\in \rS^d\F^n\setminus\{0\}$.  Then the symmetric border rank of $\cS$, denoted  as $\sbrank_{\F}\;\cS$, is $r \in \N$ if the following conditions hold
\begin{enumerate}
\item There exists a sequence $\cS_k\in\rS^d\F^n,k\in\N$ such that $\srank \cS_k=r$ for $k\in\N$ and $\lim_{k\to\infty}\cS_k=\cS$.
\item Assume that a sequence $\cS_k\in\rS^d\F^n,k\in\N$ converges to $\cS$.  Then

\noindent
$\liminf_{k\to\infty}\srank\cS_k\ge r$.
\end{enumerate}
Clearly, $\srank \cS\ge \sbrank_{\F}\;\cS$ and $\sbrank_{\F}\;\cS\ge \brank_{\F}\;\cS$.   Thus we showed
\begin{equation}\label{brankineq1}
\rank A(\cS)\le \brank_{\F}\; \cS\le \sbrank_{\F}\;\cS\le\srank \cS.
 \end{equation}

The analog of Comon's conjecture is the equality $\brank_{\F}\;\cS=\sbrank_{\F}\;\cS$.
See \cite{BGL13}.  The analog of Theorem \ref{maintheo} will be the following conjecture:
\begin{con}\label{brsymcon}  Let $d\ge 3$, $\F=\R,\C$ and $\cS\in\rS^d\F^n$.  Suppose that $\brank_{\F}\;\cS<\rank \cS$ and $\brank_{\F}\; \cS\le \rank A(\cS)+1$.
Then $\sbrank_{\F}\; \cS=\brank_{\F}\;\cS$.
\end{con}
The following theorem proves the first nontrivial case of this conjecture:
\begin{theo}\label{brank2}  Let $\cS\in\rS^d\F^n$ for $\F=\R,\C$, $d\ge 3,n\ge 2$.  Then $\brank_{\F}\;\cS=2 <\rank\cS$ if and only if
there exist two linearly independent $\x,\y\in\F^n$ and $a,b\in\F, b\ne 0$ such that 
\begin{equation}\label{brank2form}
\cS=a\otimes^d \x +b\sum_{j=0}^{d-1}(\otimes^j\x)\otimes \y\otimes (\otimes^{d-j-1}\x)
\end{equation}
In particular $\brank_{\F}\;\cS= \sbrank_{\F}\;\cS$.
\end{theo}
\subsection{Proof of Theorem \ref{brank2} }\label{Sub:prel}
\begin{lemma}\label{matrixlimlem}  Let $\F=\R,\C$ and assume that $A=[a_{i,j}]_{i\in [M],j\in[N]}\in\F^{M\times N}$, $r=\rank A$.  Suppose that the sequence
$A_k=\sum_{i=1}^q \x_{i,k} \y_{i,k}\trans, k\in\N$ satisfies the following conditions:
\begin{equation}\label{matrixlimlem1}  
\lim_{k\to\infty} A_k=A, \quad \lim_{k\to\infty}\x_{i,k}=\x_i \textrm{ for }i\in[q].
\end{equation}
\begin{enumerate}
\item  Assume that $q=r$.  Then there exists a positive integer $K$, such that for $k\ge K$ the two sets of vectors $\x_{1,k},\ldots,\x_{r,k}$ and $\y_{1,k},\ldots,\y_{r,k}$ are  linearly independent.
\item
Assume that $\x_1,\ldots,\x_q$ are linearly independent.
Then 
\begin{equation}\label{matrixlimlem2} 
\lim_{k\to\infty} \y_{i,k}=\y_i \textrm{ for } i\in[q] \quad \textrm{and } A=\sum_{i=1}^q \x_i\y_i\trans.
\end{equation}
Furthermore, $\dim\span(\y_1,\ldots,\y_q)=r$.  In particular, if $q=r$ then $\y_1,\ldots,\y_r$ are linearly independent.
\end{enumerate}
\end{lemma} 
\proof Clearly, $\rank A_k\le q$.  The first condition of \eqref{matrixlimlem1} yields that $q\ge r$.  

\noindent
 \emph{1.}   Suppose $q=r$.  Hence $\rank A_k=r$ for $k\ge K$.  Hence  the two sets of vectors $\x_{1,k},\ldots,\x_{r,k}$ and $\y_{1,k},\ldots,\y_{r,k}$ are  
linearly independent.

\noindent \emph{2.}
Complete $\x_1,\ldots,\x_q$ to a basis $\x_1,\ldots,\x_M$ in $\F^M$.
Let $\z_1,\ldots,\z_N$ be a basis in $\F^N$.  Hence $\x_i\z_j\trans, i\in[M], j\in[N]$ is a basis in $\F^{M\times N}$.  Therefore
$A_k=\sum_{i\in[M],j\in[N]}a_{ij,k}\x_i\z_j\trans$ for $k\in\N$.  The first equality of \eqref{matrixlimlem1} yields that $\lim_{k\to\infty} a_{ij,k}=a_{ij}$ for
$i\in[M],j\in[N]$.  The second equality of \eqref{matrixlimlem1} yields that $\x_{1,k},\ldots,\x_{q,k},\x_{q+1},\ldots,\x_M$ is a basis in $\F^M$ for $k\ge K$.
In what follows we assume that $k\ge K$.
Let $Q_k\in \gl(M,\F)$  be the transition matrix from the basis $[\x_1,\ldots,\x_q,\x_{q+1},\ldots,\x_M]$ to the basis 

\noindent
$[\x_{1,k},\ldots,\x_{q,k},\x_{q+1},\ldots,\x_M]$.  Clearly, $\lim_{k\to\infty} Q_k=I_M$.  Then 
\[A_k=\sum_{i\in[M],j\in[N]} b_{ij,k}\x_{i,k}\z_{j}.\]
Compare this equality with the assumption that $A_k=\sum_{i=1}^q \x_{i,k} \y_{i,k}\trans$ to deduce that $b_{ij,k}=0$ for $i>q$ and $\y_{i,k}=\sum_{j\in[N]}b_{ij,k}\z_j$
for $i\in[q]$.  Let $\tilde A_k=[a_{ij,k}], B_k=[b_{ij,k}]\in \F^{M\times N}$.  Then $B_k=Q_k\tilde A_k$.  Hence $\lim_{k\to\infty} 
B_k=\lim_{k\to\infty}Q_k A_k=\tilde A=[a_{ij}]$.  This shows \eqref{matrixlimlem2}.  
Since $\rank A=r$ it follows that $\dim\span(\y_1,\ldots,\y_q)=r$.  Thus if $q=r$ $\y_1,\ldots,\y_r$ are linearly independent.\qed 

Assume the assumptions of Definition \ref{defbordrank}.  Without loss of generality we can assume that  
\begin{equation}\label{decTk}
\cT_k=\sum_{i=1}^r \otimes_{j=1}^d \x_{i,j,k},\quad \|\x_{i,j,k}\|=1, i\in[r], j\in[d-1],   k\in\N.
\end{equation}  
By considering a subsequence of $k\in\N$ without loss of generality we can assume that 
\begin{equation}\label{convass}
\lim_{k\to\infty} \x_{i,j,k}=\x_{i,j} \quad i\in[r],j\in[d-1].
\end{equation}
(Here $\|\x\|$ is the Euclidean norm on $\F^n$.) 
\begin{lemma}\label{lembrrn}  Let $\cS\in\rS^d\F^n, d\ge 3, n\ge 2$.  Assume that $1<r=\rank A(\cS)=\brank_{\F}\;\cS<\rank \cS$.  
Let $\cT_k\in\otimes^d\F^n, k\in\N$
be a sequence of the form \eqref{decTk} satisfying \eqref{convass}.  Assume furthermore that $\lim_{k\to\infty} \cT_k=\cS$.  Then the tensors 
$\otimes_{j=1}^{d-1}\x_{1,j},\ldots,\otimes_{j=1}^{d-1}\x_{r,j}$ are linearly dependent.
\end{lemma}
\proof  Assume to the contrary that the tensors 
$\otimes_{j=1}^{d-1}\x_{1,j},\ldots,\otimes_{j=1}^{d-1}\x_{r,j}$ are linearly independent.
Lemma \ref{matrixlimlem} yields that $\lim_{k\to\infty} \x_{i,d,k}=\x_{i,d}$ for $i\in[r]$.  Hence $\cS=\sum_{i=1}^r \otimes_{j=1}^d \x_{i,j}$.
Thus $\rank \cS\le r$ which contradicts our assumptions.\qed
\begin{lemma}\label{matrlimfor}  Let $A\in\F^{M\times N}$ be a matrix of rank two.  Assume that $A_k=\a_k\b_k\trans-\c_k\d_k\trans,k\in\N$
converges to $A$.  Suppose furthermore that 
\begin{eqnarray}\label{matrlimfor1}
&&\lim_{k\to\infty} \frac{1}{\|\a_k\|}\a_k=\a,\;\lim_{k\to\infty} \frac{1}{\|\c_k\|}\c_k=\c,\; \c=\alpha\a \textrm{ for } |\alpha|=1,\\
&&\lim_{k\to\infty} \frac{1}{\|\b_k\|}\b_k=\b,\;\lim_{k\to\infty} \frac{1}{\|\d_k\|}\d_k=\d.\notag
\end{eqnarray}
Then 
\begin{equation}\label{matrlimfor1}
\b=\alpha \d, \; \lim_{k\to\infty} \|\a_k\|\|\b_k\|= \lim_{k\to\infty} \|\c_k\|\|\d_k\| =\infty, 
\; \lim_{k\to\infty} \frac{\|\c_k\|\|\d_k\|}{ \|\a_k\|\|\b_k\|}=1.
\end{equation}
Furthermore $A=\a\f\trans + \g \b\trans$, where $\span(\a,\g)=\Range A$ and $\span(\b,\f)=\Range A\trans$. 
In particular, $\g$ and $\mathbf{f}$ are limits of linear combinations of $\a_k,\c_k$ and $\b_k,\d_k$ respectively.
\end{lemma} 
\proof  Observe that 
\begin{equation}\label{Akrank2form}
A_k\trans = (\frac{1}{\|\b_k\|}\b_k)(\|\b_{k}\|\|\a_k\|)(\frac{1}{\|\a_k\|}\a_k\trans) -  (\frac{1}{\|\d_k\|}\d_k)(\|\d_{k}\|\|\c_k\|)(\frac{1}{\|\c_k\|}\c_k\trans).
\end{equation}
 Suppose first that $\span(\b)\ne \span(\d)$.  Then $\b$ and $\d$ are linearly independent.  Lemma \ref{matrixlimlem}  yields that $A\trans = a\b\a\trans + c\d\a\trans$.
Hence $\rank A=1$ which contradicts our assumptions.  As $\|\c\|=\|\d\|=1$ it follows that $\b=\beta\d$ for some scalar $\beta$ of length $1$.

We next observe that $\a\in\Range (A)$ and $\b\in\Range A\trans$.  Indeed, without loss of generality, we can assume that $\rank A_k=2$ for $k\in\N$.
Hence $\a_k\in\Range(A_k),\b_k\in\Range(A_k\trans)$.  As $\lim_{k\to\infty} A_k=A$ the assumptions \eqref{matrlimfor1} yield that 
$\a\in\Range A, \b\in\Range A\trans$.  

Assume that the sequence $\{\|\a_k\|\|\b_k\|\},k\in\N$ contains a bounded subsequence $\{n_k\},k\in\N$.  Since $\lim_{k\to\infty} A_k=A$ it follows 
the subsequence $\|\c_{n_k}\|\|\d_{n_k}\|,k\in\N$ is also bounded.  Taking convergent subsequences of the above two subsequences we deduce that
$A=\gamma \a\b\trans$.  This contradicts our assumption that $\rank A=2$.  Hence the second equality of \eqref{matrlimfor1} holds.
Rewrite \eqref{Akrank2form} as
\[A_k=\|\a_k\|\|\b_k\|((\frac{1}{\|\a_k\|}\a_k)(\frac{1}{\|\b_k\|}\b_k\trans) -
( \frac{\|\c_k\|\|\d_k\|}{ \|\a_k\|\|\b_k\|})(\frac{1}{\|\c_k\|}\c_k)(\frac{1}{\|\d_k\|}\d_k\trans)).\]
Use the assumptions that $\lim_{k\to\infty}A_k=A$, where $\rank A=2$, the facts that
$\|\a\|=\|\b\|=\|\c\|=\|\d\|=1$ and $\c=\alpha\a,\b=\beta\d$ to deduce the third and the 
the first part of \eqref{matrlimfor1}.

It is left to show that $A=\a\f\trans+\g\b\trans$.
Choose orthonormal bases $\x_1,\ldots,\x_N$ and $\y_1,\ldots,\y_N$ in $\F^M$ and $\F^N$ respectively with the following properties:
\[\x_1=\a,\; \span(\x_1,\x_2)=\Range A,\; \y_1=\b, \;\span(\y_1,\y_2)=\Range A\trans.\]
In what follows we assume that $k\gg 1$.
Choose orthonormal bases $\x_{1,k},\x_{2,k}$ and $\y_{1,k},\y_{2,k}$ in $\Range A_k$ and $\Range A_k\trans$ respectively such that
\[\x_{1,k}=\frac{1}{\|\a_k\|}\a_k, \; \y_{1,k}=\frac{1}{\|\b_k\|}\b_k,\; \lim_{k\to\infty}\x_{2,k}=\x_2,\;\lim_{k\to\infty}\y_{2,k}=\y_2.\]
Observe next that  $\{\x_{1,k},\x_{2,k},\x_3,\ldots,\x_M\}$ and $\{\y_{1,k},\y_{2,k},\y_3,\ldots,\y_N\}$ are bases in $\F^M$ and $\F^N$
which converge to bases $\x_1,\ldots,\x_N$ and $\y_1,\ldots,\y_N$ respectively.

In the bases $\{\x_{1,k},\x_{2,k},\x_3,\ldots,\x_M\}$ and $\{\y_{1,k},\y_{2,k},\y_3,\ldots,\y_N\}$ the rank one matrices $\a_k\b_k\trans, \c_k\d_k\trans$
are represented by the following block diagonal matrices: $\tilde C_k=C_k\oplus 0, \tilde D_k=D_k\oplus 0$ where
\[C_k=\left[\begin{array}{cc}a_{k}&0\\0&0\end{array}\right], \quad D_k=\left[\begin{array}{cc}b_{k}&c_k\\d_k&e_k\end{array}\right].\]
Note that $a_k=\|\a_k\|\|\b_k\|$.  Hence $\lim_{k\to\infty} a_k=\infty$.  As $\lim_{k\to\infty} A_k=A$ the arguments of the proof of Lemma \ref{matrixlimlem}
yield that $\lim_{k\to\infty} C_k-D_k=E=[e_{ij}]\in\F^{2\times 2}$.   Hence $\lim_{k\to\infty} b_k=\infty$.  Therefore $D_k=b_k\u_k\v_k\trans$,
where $\u_k\trans=(1,u_k), \v_k\trans=(1,v_k)$.  Furthermore 
\[\lim_{k\to\infty} b_ku_k=e_{21}\Rightarrow \lim_{k\to\infty} u_k=0, \quad \lim_{k\to\infty} b_kv_k=e_{12}\Rightarrow \lim_{k\to\infty} v_k=0.\]
Finally, observe that $e_{22}=\lim_{k\to\infty} b_k u_k v_k=0$.  This yields that $E=(1,0)\trans \u\trans+ \v  (1,0)$ for some choice of 
$\u,\v\in\F^2$.  Hence $A=\a\f\trans+\g\b\trans$ as claimed.

As $\rank A_k=2$ and $\lim_{k\to \infty} A_k=A$ it follows that $\g$ and $\mathbf{\f}$ are limits of linear combinations of $\a_k,\c_k$ and $\b_k,\d_k$ respectively.
\qed

\textbf{Proof of Theorem \ref{brank2}.}   Assume first that $\cS$ is of the form \eqref{brank2form} where $b\ne 0$. 
Without loss of generality we can assume that $b=1$.
Clearly, $\Range A(\cS)=\span(\x,\y)$  Hence $\rank A(\cS)=2$.
Let $\cT(\epsilon):=a\otimes^d \x+\frac{1}{\epsilon}(\otimes^d(\x+\epsilon\y)-\otimes^d\x)$ for $\epsilon\ne 0$.  
Then $\rank\cT(\epsilon)=2$ for $\epsilon^{-1}\ne a$.  Clearly, $\lim_{\epsilon\to 0}\cT(\epsilon)=\cS$.  Hence $\brank_{\F}\;\cS=2$.  As $\cT(\epsilon)\in\rS^d\F^n$ 
it follows that $\sbrank_{\F}\;\cS=2$.
We claim that $\rank\cS>2$.  We can assume without loss of generality that $n=2$ and $\x=\e_1=(1,0)\trans,\y=\e_2=(0,1)\trans$.
 Assume first that $d=3$.   So $\cS=[s_{i,j,k}]$ where
\[s_{1,1,1}=a,\; s_{1,1,2}=s_{1,2,1}=s_{2,1,1}=1,\; s_{1,2,2}=s_{2,1,2}=s_{2,2,1}=s_{2,2,2}=0.\]
Let 
\[F=[s_{i,j,1}]_{i,j\in[2]}=\left[\begin{array}{rr}a&1\\1&0\end{array}\right], \quad G=[s_{i,j,2}]_{i,j\in[2]}=\left[\begin{array}{rr}1&0\\0&0\end{array}\right].\]
Then $\rank\cS=2$ if and only if the matrix $GF^{-1}$ is diagonalizable, see e.g. \cite{Fri13}.
Clearly, $GF^{-1}=\left[\begin{array}{rr}0&1\\0&0\end{array}\right]$ is not diagonalizable.  Hence $\rank \cS>2$.  It is easy to show straightforward that $\rank\cS=3$.

Assume now that $d>3$.   Let $\phi:\F^2\to\F$ be the linear functional such that $\phi(\e_1)=\phi(\e_2)=1$.
Consider the following map $\psi:(\F^2)^d\to \otimes^3\F^2$:
$\psi((\u_1,\ldots,\u_d))=(\prod_{j=4}^d \phi(\u_j))\otimes_{i=1}^3\u_j$.
Clearly, $\psi$ is a multilinear map.  The universal lifting property of the tensor product yields that $\psi$ lifts to the linear map $\Psi: \otimes^d \F^2\to \otimes^3\F^2$ such that
\[\Psi(\otimes_{j=1}^d\u_i)=(\prod_{j=4}^d \phi(\u_j))\otimes_{i=1}^3\u_j.\]
Observe that a rank one tensor is mapped to either rank one tensor or zero tensor.  Clearly, the image of a symmetric rank
one tensor is a symmetric tensor of at most rank one.  Hence
$\Psi:\rS^d\F^2\to\rS^3\F^2$.  Assume that $\cS$ of the form  \eqref{brank2form}, where $\x=\e_1$ and $\y=\e_2$.  Then
\[\Psi(\cS)=(a+(d-3)b)\x\otimes\x\otimes\x+b(\y\otimes \x\otimes \x+\x\otimes\y\otimes \x +\x\otimes\x \otimes \y).\]
Assume to the contrary that $\rank\cS=2$. Then $\rank \Psi(\cS)\le 2$.
This contradicts our proof that $\rank \Psi(\cS)=3$.  Hence $\rank \cS\ge 3$.

Assume now that $\cS\in\rS^d\F^n$ and $2=\rank A(\cS)=\brank_{\F}\;\cS<\rank \cS$.
Let $\cT_k\in\otimes^d\F^n$ be a sequence of tensors of rank two converging to $\cS$.  So $\cT_k=\otimes_{j=1}^d \x_{j,k}-\otimes_{j=1}^d\y_{j,k}$.
Since $\rank A(\cS)=2$ we can assume without loss of generality: First, $\x_{j,k}$ and $\y_{j,k}$ are linearly independent for $j\in[d],k\in\N$.   Second,
\[\lim_{k\to\infty} \frac{1}{\|\x_{j,k}\|} \x_{j,k}=\x_{j},\;\lim_{k\to\infty} \frac{1}{\|\y_{j,k}\|} \y_{j,k}=\y_j \textrm{ for }j\in[d].\]
Lemma \ref{lembrrn} yields that $\otimes_{j=1}^{d-1} \x_j$ and $\otimes_{j=1}^{d-1} \y_j$ are linearly dependent.  
Hence $\span(\x_j)=\span(\y_j)$ for $j\in[d-1]$.
Lemma \ref{matrlimfor} yields that $\x_d$ and $\y_d$ linearly dependent.  So $\span(\x_d)=\span(\y_d)$.   
Apply Lemma \ref{matrlimfor} to $A=A(\cS)\trans, A_k=A(\cT_k)\trans,k\in \N$.
It then follows that  $\cS=(\otimes_{j=1}^{d-1}\x_j)\otimes \z +\cF\otimes \x_d$ 
for some $\cF\in \otimes^{d-1}\F^n$.  
Furthermore $\span(\x_d,\z)=\Range A(\cS)$.   As $\x_d$ and $\z$ are linearly independent and $\cS$ symmetric it follow that $\otimes_{j=1}^{d-1},\cF\in
\rS^{d-1}\F^n$.  Hence $\span(\x_1)=\cdots=\span(\x_{d-1})=\span(\x)$.   Thus $\otimes_{j=1}^{d-1}\x_j=t\otimes^{d-1}\x$.  
By considering the unfolding of $\cS$ in another mode we deduce that $\span(\x_d)=\span(\x)$.  Observe
next that $\rank \cF>1$.  Otherwise $\rank \cS\le 2$ which contradicts our assumptions.  Lemma \ref{matrlimfor} yields that  $\cF$ is a limit of linear
combinations of $\otimes_{j=1}^{d-1} \x_{j,k}$ and  $\otimes_{j=1}^{d-1} \y_{j,k}$.
Hence $\brank_{\F} \cF\le 2$.  As $\rank \cF>1$ it follows that $\brank_{\F}\; \cF= 2$. 
In summary we showed:
\begin{eqnarray}\label{tempforS}
&&\cS=\otimes^{d-1}\x\otimes \z +\cF\otimes \x,\\ 
&&\span(\x,\z)=\Range A(\cS),\; \cF\in\rS^{d-1}\F^n,\;\brank_{\F}\;\cF=2.\notag
\end{eqnarray} 

We now prove the following claim:  Assume that $\cS\in\rS^d\F^n, \rank A(\cS)=2<\rank \cS$.  Suppose furthermore $\cS$ is a limit of linear
linear combinations of $\otimes_{j=1}^d \x_{j,k},\otimes_{j=1}^d \y_{j,k}$, where the the following limit exist and satisfy:
\begin{equation}\label{limitconxy}
\lim_{k\to\infty} \frac{1}{\|\x_{j,k}\|} \x_{j,k},\lim_{k\to\infty} \frac{1}{\|\y_{j,k}\|} \y_{j,k}\in \span(\x)\textrm{ for }j\in[d].
\end{equation}
Then \eqref{brank2form} holds.

We prove the claim by induction on $d$.  Assume first that $d=3$.  Observe first that two dimensional
subspace $\Range A(\cS)\trans=\span(\a,\c)$, as given by Lemma \ref{matrlimfor}, is in $\rS^2\F^n$. i.e. the space of symmetric matrices.
Lemma \ref{matrlimfor} yields that $\cF$ is a limit of linear combinations of
$\x_{1,k}\otimes \x_{2,k}$ and $\y_{1,k}\otimes \y_{2,k}$. As 
$\lim_{k\to\infty}\frac{1}{\|\x_{1,k}\|\|\x_{2,k}\|}\x_{1,k}\otimes \x_{2,k}=t\x\otimes \x$ it follows that $\Range A(\cS)\trans$
contains rank one matrix $\x\otimes\x$.  Lemma  \ref{matrlimfor} yields that  $\Range A(\cS)\trans$ contains a rank two matrix of the form $\x\otimes \mathbf{f}+
\g\otimes\x$.  Since this matrix is symmetric it is of the form $c\x\otimes\x +\x\otimes \u+\u\otimes \x$ for some scalar $c$ and $\u\in\F^n$ which is linearly
independent of $\x$.  Hence $\cF=d\x\otimes\x +  \x\otimes\v+\v\otimes \x$ for $\v=d\u, d\ne 0$.  As rank $A(\cS)=2$ it follows that $\span(\x,\v)=\span(\x,\z)$.
Therefore we showed that
\[\cS=a\x\otimes\x\otimes\x +b\x\otimes\x\otimes \v+c(\x\otimes \v+\v\otimes \x)\otimes \x.\] 
Interchange the last two factors in $\cS$ to deduce that $0=(b-c)(\x\otimes \x\otimes \v-\x\otimes\v\otimes \x)$.  
Hence $b=c$ and \eqref{brank2form} holds for $d=3$.

Assume now that \eqref{brank2form} holds for $d=p$ and suppose that $d=p+1$.  Consider \eqref{tempforS}.  
Suppose first that $2<\rank\cF$.
Then the induction hypothesis applies to $\cF$.  Hence $\cF$ is of the form \eqref{brank2form} and
\[\cS=\otimes^p\x\otimes \z+(a\otimes^p\x+b \sum_{j=0}^{p-1}\otimes^j\x\otimes\y\otimes (\otimes^{p-1-j}\x))\otimes \x, \;b\ne 0.\]
Note that $\Range A(\cS)=\span(\x,\y)=\span(\x,\z)$.  Hence
\[\cS=a'\otimes^{p+1}\x +c \otimes^p\x\otimes \y +b \sum_{j=0}^{p-1}\otimes^j\x\otimes\y\otimes (\otimes^{p-j}\x).\]
Interchange the last two factors in $\cS$ to deduce that $(c-b)\otimes^{p-1}\x(\x\otimes\y -\y\otimes \x)=0$.  Hence $b=c$ and $\cS$ is of the  form
 \eqref{brank2form}.

It is left to consider the case where $\cF$ is a symmetric tensor of rank two.  
So $\rank \cF=2$.  Hence $\rank A(\cF)=2$.  Theorem \ref{eqcase}   yields that $\cF=s  \otimes^p \u +t \otimes ^d\v$, were $s,t=\pm 1$, and this decomposition is unique.
(The $\pm $ are needed if $\F=\R$ and $p$ is even.)  Clearly, $\span(\u,\v)=\Range A(\cS)$.  It is enough to assume that $n=2$.
Recall that $A(\cF)$ is a limit of a linear combinations of two rank one matrices:
$\x_{1,k}(\otimes_{j=2}^{p}\x_{j,k})\trans, \y_{1,k}(\otimes_{j=2}^{p}\y_{j,k})\trans, k\in\N$.  The assumption \eqref{limitconxy} implies that
we can use Lemma \ref{matrlimfor}.  Hence $\cF=\x\otimes \cG+\g\otimes (\otimes^{p-1}\x)$.  Therefore $\otimes^{p-1}\x\in\span(\otimes^{p-1}\u,\otimes^{p-1}\v)$.
We claim that this possible if and only if either $\span(\x)=\span(\u)$ or $\span(\x)=\span(\v)$.    Suppose to the contrary that  $\span(\x)\ne\span(\u)$ or $\span(\x)\ne\span(\v)$.
So $\otimes^{p-1}\x=s\otimes^{p-1}\u+t\otimes^{p-1}\v$.  Clearly $st\ne 0$.
Let $\phi:\F^2\to \F$ be a nonzero linear functional such that $\phi(\x)=0$.  Let $\Psi:\otimes^{p-1}\F^2\to \otimes^{p-2}\F^2$ be the linear mapping
of the form given above: $\Psi(\otimes_{j=1}^{p-1}\w_j)=\phi(\w_1)\otimes_{j=2}^{p-1}\w_j$.  So 
\[0=\Psi(\otimes^{p-1}\x)=s\phi(\u)\otimes^{p-2}\u+t\phi(\v)\otimes^{p-2}\v.\]
This is impossible since $\u$ and $\v$ are linearly independent.  Thus we can assume that $\span(\v)=\span(\x)$,  $\cT=s\otimes^p\u+t'\otimes^p\x$
and $\z\in\span(\u,\x)$.  Thus
\[\cS=a\otimes^{p+1}\x+b\otimes^p\x\otimes \u+s\otimes^p\u\otimes \x.\]
As $\rank \cS=3$ it follows that $\cS':=\cS-a\otimes^{p+1}\x$ is a symmetric tensor of rank two.  Theorem \ref{eqcase} claims that the decomposition $\cS'=
b\otimes^p\x\otimes \u+s\otimes^p\u\otimes \x$ is unique and $\otimes^p\x\otimes \u,\otimes^p\u\otimes \x$ are symmetric tensors.  So $\span(\u)=\span(\x)$
which contradicts our assumption that $\u$ and $\x$ are linearly independent.\qed
\subsection{Approximation of symmetric tensors}
Define on $\otimes^d\F^n$ the standard inner product:
\[\an{\cP,\cQ}:=\sum_{i_j\in[n],j\in[d]}p_{i_1,\ldots,i_d}\overline{q_{i_1,\ldots,i_d}},\quad \cP=[p_{i_1,\ldots,i_d}],\cQ=[q_{i_1,\ldots,i_p}]\in\otimes^d\F^n.\]
Observe that $\an{\otimes_{j=1}^d\x_j,\otimes_{j=1}^d\y_j}=\prod_{j=1}^d \an{\x_j,\y_j}$.
Assume that $k\in[1,d-1]$.  Denote by $\Gr(k,\F^n)$ the Grasmannian manifold of $k$-dimensional subspace in $\F^n$.  Let 
\[\Gr(k,d,\F^n):=\{\otimes_{j=1}^d \U_j, \quad \U_j\in\Gr(k,\F^n), j\in[d]\}.\]
For a given $\otimes_{j=1}^d \U_j$ denote by $P_{\otimes_{j=1}^d \U_j}:\otimes^d \F^n \to \otimes_{j=1}^d \U_j$ the orthogonal projection
of $\otimes^d \F^n$ on $\otimes_{j=1}^d \U_j$. A best $k$-approximation of $\cT\in\otimes^d\F^n $ is each tensor $\cT^\star$ satisfying
\[\min_{\otimes_{j=1}^d \U_j\in\Gr(k,d,\F^n)}\|\cT-P_{\otimes_{j=1}^d \U_j}(\cT)\|=\|\cT-\cT^\star\|, \; \cT^\star=P_{\otimes_{j=1}^d \U_j^\star}(\cT),
\otimes_{j=1}^d \U_j^\star\in\Gr(k,d,\F^n).\]
See \cite{FS13, FT15}.  The results of \cite{FS14} yield that $\cT^\star$ is unique for $\cT$ outside of a semi-algebraic set of dimension less than the real dimension of
$\otimes^d\F^n$.
The analog of Comon's conjecture is:
\begin{con}\label{Comapcon}  Let $n-1,d-1\in\N, k\in[d-1]$ and $\cS\in \rS^d\F^n$.  Then a best $k$-approximation of $\cS$ can be chosen to be a symmetric tensor.
\end{con} 
This conjecture is known to hold in the following cases:
For $d=2$ it is a consequence of Singular Value Decomposition.  For $k=1$ and $d>2$ it follows from Banach's theorem \cite{Ban38}.
See \cite{Fri13} for $\F=\R$.  Similar arguments combined with Banach's theorem yield the case $\F=\C$.
It is shown in \cite{FS13} that for $\F=\R$ there is a semi-algebraic set in $\rS^d\R^n$ of dimension $\dim \rS^d\R^n$ for which
the conjecture holds.\\

\emph{Acknowledgment}: I thank Lek-Heng Lim for bringing my attention to Comon's conjecture, to  Joseph M. Landsberg
and the anonymous referees for their useful comments.

\bibliographystyle{plain}

\end{document}